\documentclass[11pt]{article}
\usepackage{geometry, amsmath, amsfonts, amssymb,pdfsync} 
\usepackage{graphicx, graphics, psfrag}    
\usepackage{latexsym}
\usepackage{pdfsync}
\usepackage{mathrsfs,color}           
\newtheorem{theorem}{Theorem}
\newtheorem{proposition}{Proposition}[section]
\newtheorem{lemma}[proposition]{Lemma}
\newtheorem{corollary}[proposition]{Corollary}
\newtheorem{definition}[proposition]{Definition}

\newtheorem{remark}[proposition]{Remark}

\numberwithin{equation}{section}
\newcommand{\bpf} {\noindent{\sc Proof} : }
\newcommand{\epf}  {\hfill$\diamondsuit$\vspace{.5cm}}
\newcommand{\E} {\mathbb{E}}
\newcommand{\EE} {\mathcal{E}}
\newcommand{\R} {\mathbb{R}}
\renewcommand{\P} {\mathbb{P}}

\newcommand{\Q} {\mathbb{Q}}
\newcommand{\N}{\mathbb{N}}
\newcommand{\F} {\mathcal{F}}
\newcommand{\G} {\mathcal{G}}
\newcommand {\eps} {\varepsilon}

\title{A path--valued Markov process \\ indexed by the ancestral mass\thanks{We acknowledge support by MAEDI/MENESR and DAAD through the PROCOPE programme. A.W. was supported in part by DFG SPP 1590 {\it Probabilistic Structures in Evolution.}}}
\author{Etienne Pardoux\footnote{Aix-Marseille Universit\'e, CNRS, Centrale Marseille, I2M, UMR 7373 13453 Marseille, France \newline
\tt{e-mail: etienne.pardoux@univ-amu.fr} \newline \tt{URL:  http://www.i2m.univ-amu.fr/$\sim$pardoux/ }} 
\and Anton Wakolbinger\footnote{Goethe-Universit\"at, Institut f\"ur Mathematik, Robert-Mayer-Str. 10, 60325 Frankfurt, Germany. \newline
\tt{e-mail: }\tt {wakolbinger@math.uni-frankfurt.de} \newline \tt{URL: http://www.math.}\tt{uni-frankfurt.de/}\tt{$\sim$ismi/wakolbinger/pers.html}}}
\date{}
\begin{document}
\maketitle
\begin{abstract}
A family of Feller branching diffusions $Z^x$, $x \ge 0$, with nonlinear drift and initial value $x$ can, with a suitable coupling over the {\em ancestral masses} $x$, be viewed as a path-valued process indexed by $x$. For a coupling due to Dawson and Li, which in case of a linear drift describes the corresponding Feller branching diffusion, and in our case makes the path-valued  process Markovian, we find an SDE solved by $Z$, which is driven by a random point measure on excursion space. In this way we are able to identify the infinitesimal generator of the path-valued process. We also establish  path properties of $x\mapsto Z^x$ using various couplings of $Z$ with classical Feller branching diffusions.
\end{abstract}
{\small {\bf Subject classification:} 60J25, 60J75, 60J80, 92D25.} \\
{\small {\bf Keywords and phrases:} Feller branching process with interaction, point process of excursions, Girsanov transform.}
\section{Introduction}
Consider the
SDE
\begin{equation}\label{eqDL}
Z^x_t=x+\int_0^tf(Z^x_s)ds+2\int_0^t\int_0^{Z^x_s}W(ds,d\xi), \quad t\ge 0, \, x \ge 0
\end{equation}
where
$W(\cdot,\cdot)$ denotes a two--dimensional white noise,
i.e. a generalized zero mean Gaussian random field on $\R^2_+$, whose covariance operator is the identity operator
on $L^2(\R_+^2)$. Our assumptions on the function $f: \mathbb R_+ \to \mathbb R$ will be specified below.

For any fixed $x>0$, the solution of \eqref{eqDL} has the same law as the solution of  the simpler SDE 
$$Z^x_t=x+\int_0^tf(Z^x_s)ds+2\int_0^t\sqrt{Z^x_s}dB_s,$$
 where $B$ is a standard scalar Brownian motion. However, the formulation \eqref{eqDL}, which follows Dawson and Li \cite{DL}, is the simplest and most natural way to specify 
 a coupling which in the case of a linear drift renders the branching property of $Z$, i.e. the independence of $Z^x$ and $Z^{x+y}-Z^x$ for all $x,y > 0$.
 For a quadratic $f$, this same coupling is discussed also in \cite{PW2} Sec. 2, but there without appealing to the representation \eqref{eqDL}.

When $f$ is linear, i.e. $f(z)=\theta z$, the solution $Z^x_t$, $t\ge 0$,  is a continuous
state branching process with ancestral mass $x$.  This process models the evolution of the (continuum scaling limit of the) size of a population when the reproduction dynamics of the various individuals are mutually independent. 

In the terminology of a (pre-limiting) individual-based description, the possible nonlinearities of $f$ model an impact of the current population size on the individual reproduction dynamics, and  in this way go along with an interaction in the individuals' reproductive behavior.
 If the interaction is of 
the type of competition for rare resources, then an increase in the population size decreases the individual birth rate and/or increases the death rate.
This means that $f(z)/z$ should be decreasing, and $f(z)$ should be negative  for $z$ large enough. On the other hand, specifically for moderate values of $z$,
$f(z)/z$ might be increasing. This is the case in the presence of the so--called Allee effect, where there is a negative growth rate for small population sizes $z$ and a positive growth rate for larger population sizes (as long as the population size does not exceed a certain carrying capacity).

In previous publications, we obtained several results on the solution of
equation \eqref{eqDL} (for $f$ of the form $f(z)=\theta z-\gamma z^2$ or for more general $f$). In particular we discussed

(i)   its approximation by finite population models (\cite{LPW}, \cite{BP}),

(ii)  the extension of the second Ray--Knight theorem and a description of the forest of genealogical trees of the population whose total size follows \eqref{eqDL} (\cite{LPW}, \cite{BP}, \cite{PW2}),

(iii)  the effect of competition on the asymptotic extinction time and total mass of the forest of trees for large population size (\cite{LP}).

In this paper, we study the solution of \eqref{eqDL} as a path-valued process indexed by the mass $x$ of the ancestral population. 

In the case where $f$ is linear,
much is known about the $\R_+$--valued process $\{Z^x_t,\ x>0\}$ for any $t>0$ fixed, as well
as about the $C(\R_+;\R_+)$--valued process $\{Z^x_\cdot,\ x>0\}$. Those have independent increments, which reflects the independence 
of the  progenies of various ancestors in a branching process.  Moreover, for any fixed $t>0$, $x\mapsto Z^x_t$ is an increasing process which
has a.s. finitely many jumps on any finite time interval, and is constant between its jumps. 
On the other hand, the path-valued process $x\mapsto Z^x_\cdot$ has infinitely many jumps on any  interval of positive Lebesgue measure. Also, for linear $f$ explicit formulas for the law of the random variables $Z_t^x$ are available.

In the presence of a nonlinear drift $f$, 
the situation is more complicated, and our paper contributes to its investigation. Our first step is to show that $x\mapsto Z^x_t$ is, for fixed $t$,  again 
an increasing process which increases only by jumps, whose number is finite on any compact interval. Now the process $x\mapsto Z^x_t$ 
in general does no longer have independent increments. We strongly suspect that (for fixed $t$) it is not a Markov process. However, the path-valued process $x\mapsto Z^x_\cdot$ {\em is} Markovian, and it is the objective of this paper to write  an SDE driven by a random point measure for the process  $x\mapsto Z^x_\cdot$, and to identify its infinitesimal generator, by writing a martingale problem formulation of the path--valued SDE which it solves, see Theorem \ref{th:path-val-SDE} and Corollary \ref{cor:martpb} below.

Let us now specify our standing assumptions on the nonlinear function $f$.
We assume that $f\in C(\R_+;\R)$, $f(0)=0$,  and $f$ satisfies in addition the  following three assumptions
\begin{align}
f&(a+b)-f(a)\le\theta b,\quad \text{for some }\theta\ge0 \text{ and all }a,b>0\, ;\label{1}\\
\begin{split}
f \text{  is }\frac{1}{2}\text{--H\"older} & \text{ continuous, i.e. for all $M > 0$ there exists a } C_M < \infty \text{ such that } 
\\ 
|f &(a+b)-f(a)| \le C_M\sqrt b\, \text{ for all } a \in [0,M] \text { and } b \in [0,1]\, ;  \label{2}
\end{split}
\\
\int_1^\infty&\exp\left(-\frac{1}{2}\int_1^u\frac{f(r)}{r}dr\right)du=+\infty.\label{3}
\end{align}
The first assumption is crucial for equation \eqref{eqDL} to be well--posed, see \cite{BP}. It also implies that $f(z)\le \theta z$ for all $z>0$, which 
will be important in the next section. 
The second assumption implies that $b^{-1/2}[f(a+b)-f(a)]$ remains bounded while $b\to0$, which will be essential for our Girsanov transformations below. Finally, it is shown in \cite{BP} that the third assumption implies that, for all $x >0$, the random path $Z^x_\cdot$ hits zero in finite time a.s., which will be important in many of our arguments below. For technical as well as for conceptual reasons we want to study models of populations which go extinct in finite time. This is why we assume that $f(0)=0$, and not just that $f(0)\ge0$. 
In particular, we do not consider populations with immigration (except as an auxiliary construction in some proofs below).
 
Note that a sufficient condition for \eqref{3} to hold is that there exists $z_0>0$ such that $f(z)\le2$, for all $z\ge z_0$. Clearly, a wide variety of  functions $f$ satisfy our assumptions.
\bigskip

The paper is organized as follows. In section \ref{se1}, we establish the basic properties of the solution of \eqref{eqDL}, recalling in particular the existence and uniqueness result from \cite{DL}. 

Section \ref{se2} is devoted to comparison with a supercritical Feller diffusion $Y$, with supercriticality parameter $\theta$, the same real number which appears in the assumption \eqref{1}. We first prove a basic and easy comparison theorem between $Z$ and $Y$. Next we construct another coupling of the two processes, for
which a much stronger comparison holds. This permits us to deduce that $Z$ increases only where $Y$ increases, in particular $x\mapsto Z^x |_{[\delta, \infty)}$ is constant between its jumps for all $\delta > 0$.

Section \ref{sec:3} is devoted to establishing a path--valued SDE satisfied by $\{Z^x_\cdot,\  x>0\}$, and deducing the exact form of the generator of that Markov process. Here again we shall consider a pair $(Y,Z)$. However the process $Y$ will then be a critical Feller diffusion, and in this case there will be no comparison between $Z$ and $Y$. Instead, we shall  exploit Girsanov's theorem  and write the Radon--Nikodym derivative of the law of $Z$ with respect to that of $Y$. This will be our key ingredient for the identification of the generator of the process $Z^x_\cdot$.

\section{Basic results}\label{se1}
It follows from Theorem 2.1 in \cite{DL} that for any $x>0$, \eqref{eqDL} has a unique continuous non--negative solution. Note in particular that the assumptions of that theorem are satisfied here, since we can decompose $f(a)=\theta a+[f(a)-\theta a]$, where $a\mapsto \theta a$ is  Lipschitz, while 
$a\mapsto f(a)-\theta a$ is continuous and non increasing. 

Since we are interested in the two parameter process $\{Z^x_t,\ t\ge0,\ x>0\}$, we need to make sure that we can choose an appropriate version.
\begin{lemma}\label{cont-proba}
The mapping $\xi\mapsto Z^\xi_\cdot$ is continuous in probability.
\end{lemma}
\bpf
Let $x,y>0$. Theorem~2.2 in \cite{DL} implies that
  \begin{equation}\label{compsimple}
  \P(Z^{x+y}_s-Z^x_s\ge 0,\,  \forall s \ge 0) = 1\,. 
  \end{equation} 
  Consequently,
\begin{equation}\label{Ztinc}
Z^{x+y}_t-Z^x_t=y+\int_0^t[f(Z^{x+y}_s)-f(Z^x_s)]ds+2\int_0^t\int_{Z^x_s}^{Z^{x+y}_s}W(ds,d\xi).
\end{equation}
 Taking the expectation in this identity, and exploiting \eqref{1} and Gronwall's Lemma, we infer that
 \begin{equation}\label{EZtinc}
 \E[Z^{x+y}_t-Z^x_t]\le y\exp(\theta t).
 \end{equation}
 Let $M_t$ denote the last term on the right of \eqref{Ztinc}. The quadratic variation of this martingale is given by
 $$\langle M,M\rangle_t=4\int_0^t (Z^{x+y}_s-Z^x_s)ds.$$
 A consequence of \eqref{Ztinc} and \eqref{1} is that for any $t>0$~:
 \begin{equation}\label{supZinc}
 \sup_{0\le s\le t} (Z^{x+y}_s-Z^x_s)\le y+\theta\int_0^t (Z^{x+y}_s-Z^x_s)ds
 +\sup_{0\le s\le t}|M_s|.
 \end{equation}
 Taking the expectation in \eqref{supZinc}, we deduce from 
 the Burkholder--Davis--Gundy and Schwartz inequalities and \eqref{EZtinc} that there exists a constant $c>0$ such that
 $$\E\left[\sup_{0\le s\le t} (Z^{x+y}_s-Z^x_s)\right]\le y(1+e^{\theta t})+c\sqrt{y\theta^{-1}e^{\theta t}},$$
 from which the result follows, again in view of \eqref{compsimple}.
   \epf
 
  \begin{definition}
  We will denote by $E$ the space of continuous functions $u$ from $[0,+\infty)$ into itself, which are such that whenever
 $\zeta(u):=\inf\{t>0,\ u(t)=0\}$ is finite, then $u(t)=0$, for any $t\ge\zeta(u)$.  We equip $E$ with the topology of uniform convergence on compacts.  
 \end{definition}

The following result is similar to Theorem 3.6 in \cite{DL}.
\begin{lemma}\label{cadlag}
There exists a version of the mapping $\xi\mapsto Z^\xi_\cdot$ which is a.s. increasing and c\`adl\`ag with values in $E$.
\end{lemma}
\bpf
Let $\{\tilde{Z}^\xi,\ \xi>0\}$ denote an arbitrary collection of processes, such that $\tilde{Z}^\xi$ solves the SDE
\eqref{eqDL} for any $\xi>0$.
The fact that $\xi\mapsto \tilde{Z}^\xi_\cdot$ is a.s. increasing from $\Q_+$ into $C([0,+\infty))$  follows from \eqref{compsimple}.
Now  for any $x>0$ we define
\[ Z^x=\lim_{\xi_n\downarrow x,\ \xi_n\in\Q_+}\tilde{Z}^{\xi_n}.\]
By monotonicity, the sequence converges a.s., and it follows from Lemma \ref{cont-proba} that, for any $x>0$, 
$Z^x=\tilde{Z}^x$ a.s., hence $Z^x$ solves \eqref{eqDL}. The result follows.
 \epf
 
 As we will see below, the mapping $\xi \mapsto Z^\xi_\cdot$ does have discontinuities with positive probability.

\section{Connection with a supercritical Feller diffusion}\label{se2}
In this section,
$\{Y^x_t,\, t\ge0,\,  x>0\}$ stands for a Feller branching diffusion with supercriticality parameter $\theta$, starting from an ancestral mass $x>0$. More precisely, for a given space--time white noise $W$, and  $\theta>0$  being the parameter that enters condition \eqref{1} on~$f$, we write $Y^x$ for the solution of 
\begin{equation}\label{eq:Y}
Y^x_t=x+\theta\int_0^tY^x_sds+2\int_0^t\int_0^{Y^x_s}W(ds,d\xi).
\end{equation}
Let $Z^x$ be the solution of \eqref{eqDL}, with $f$ satisfying conditions \eqref{1}, \eqref{2} and \eqref{3}.
The two equations \eqref{eq:Y} and \eqref{eqDL} with the same $W$ describe one possible coupling of the two random fields $\{Y^x_t,\, t\ge0, \,x>0\}$ and $\{Z^x_t,\, t\ge0,\,  x>0\}$. 
\begin{proposition}
For each $x>0$,
$$\P(Z^x_t\le Y^x_t,\, \forall t\ge0)=1.$$
\end{proposition}
\bpf
Since  $f(x)\le\theta x$, this is immediate from the comparison theorem  (Theorem 2.2) in \cite{DL}.
\epf

We now construct yet another coupling which will allow to derive distributional properties of $Z$ that are required in the sequel.
For each $t>0$, $x>0$, let
\begin{align*}
D_t&=\{\xi>0;\  Y^\xi_t>Y^{\xi-}_t\},\ \text{and}\\
A^x_t(Z)&=\cup_{\xi\le x,\ \xi\in D_t}(Y^{\xi-}_t,Y^{\xi-}_t+Z^\xi_t-Z^{\xi-}_t].
\end{align*}
Note that the random set $A^x_t$ depends upon the copy of $Z$, in particular upon the chosen coupling of $Y$ and $Z$.
Note also that the Lebesgue measure of the set $A^x_t(Z)$ equals $Z^x_t$.
\begin{figure}[h]
\begin{center}
\psfrag{a}{$\tilde Z^{\xi-}_t$}
\psfrag{c}{$Y^{\xi-}_t+ \tilde Z^{\xi}_t- \tilde Z^{\xi-}_t$}
\psfrag{b}{$Y^{\xi-}_t$}
\psfrag{d}{$Y^{\xi}_t$}
\psfrag{x}{$\xi$}
\psfrag{t}{$t$}
\psfrag{r}{$R^\xi$}
\includegraphics[width=7cm, height=9cm]{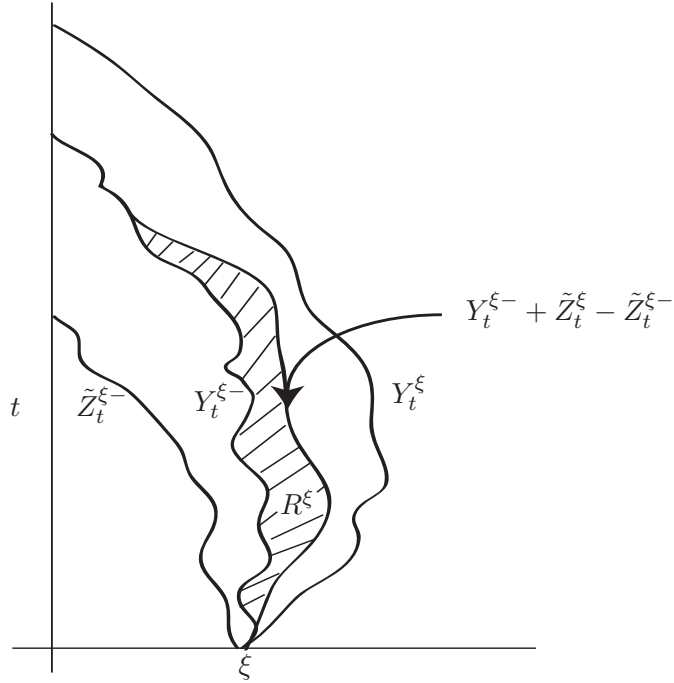}
\caption{The region $R^\xi$ of the noise $W$ that drives $\tilde Z^\xi -\tilde Z^{\xi-}$ (which is shaded in the picture) is contained in the region of the noise that drives $Y^\xi-Y^{\xi-}$ (which is the one between $Y^{\xi-}$ and  $Y^{\xi}$). In particular, $R^\xi$ does not intersect the region to the left of $Y^{\xi-}$, which would be the case if $\tilde Z$ were replaced by $Z$. }
\end{center}
\end{figure}
We have the
\begin{theorem}\label{def-coupling}
There exists a random field $\{\tilde Z^x_t,\ x>0,t\ge0\}$ such that $t\mapsto\tilde Z^x_t$ is continuous, $x\mapsto\tilde Z^x_t$ is right--continuous,
$\{\tilde Z^x_t,\ x>0,t\ge0\}$ has the same law as $\{Z^x_t,\ x>0,t\ge0\}$ (the solution of \eqref{eqDL}), $\{\tilde Z^x_t,\ x>0,t\ge0\}$
solves the SDE
\begin{equation}\label{eq:DLcoupled}
\tilde Z^x_t=x+\int_0^t f(\tilde Z^x_s)ds+2\int_0^t\int_{A^x_s(\tilde Z)}W(ds,d\xi),
\end{equation}
and moreover for all $x,y>0$,
\begin{equation}\label{comparStrong}
\P(\tilde Z^{x+y}_t-\tilde Z^x_t\le Y^{x+y}_t-Y^x_t,\, \forall t\ge0)=1.
\end{equation}
\end{theorem}
\bpf
For a solution $\tilde Z$ of \eqref{eq:DLcoupled}, the equality in law between $\{\tilde Z^x_t,\ x>0,t\ge0\}$ and $\{Z^x_t,\ x>0,t\ge0\}$ follows from the fact that the Lebesgue measure
of $A^x_t(\tilde Z)$ equals $\tilde Z^x_t$. We now construct a solution of \eqref{eq:DLcoupled}.

For each $k,n\ge1$, let $x^k_n:=2^{-n}k$.
For each $n\ge1$, we now define  $\{Z^{n,x}_t,\, t\ge0\}$.
For $0<x\le x^1_n$, we require that $\{Z^{n,x}_t,\ t\ge0\}$ solves
\begin{equation*}\label{eq:n,1}
Z^{n,x}_t=x+\int_0^tf(Z^{n,x}_s)ds+2\int_0^t\int_0^{Z^{n,x}_s}W(ds,d\xi).
\end{equation*}
And for  $k\ge2$, we define recursively $\{Z^{n,x}_t,\ t\ge0\}$ for $x^{k-1}_n<x\le x^k_n$ as the solution of
\begin{equation*}\label{eq:n,k}
Z^{n,x}_t-Z^{n,{x^{k-1}_n}}_t=x-x^{k-1}_n+\int_0^t\left[f(Z^{n,x}_s)-f(Z^{n,{x^{k-1}_n}}_s)\right]ds
+2\int_0^t\int_{Y_s^{x^{k-1}_n}}^{Y_s^{{x^{k-1}_n}}+Z^{n,x}_s-Z^{n,{x^{k-1}_n}}_s}\!\!\!\!\!\!\!\!\!\!\! W(ds,d\xi).
\end{equation*}
From \eqref{1} together with \eqref{eq:Y} and Theorem 2.2 in \cite{DL} it follows that  for all $k\ge1$ and $x^{k-1}_n<x\le x^k_n$,
\begin{align}\label{compar}
Z^{n,x}_t-Z^{n,{x^{k-1}_n}}_t\le Y^x_t-Y^{x^{k-1}_n}_t \mbox{ a.s. for all } t\ge0\, .
\end{align}
Moreover, 
 the law of $\{Z^{n,x}_t,\, x>0, t\ge0\}$ is the same as the law of
 $\{Z^x_t,\, x>0,\, t\ge0\}$, the solution of \eqref{eqDL}. 
 
 Recall that for each $t>0$, $x\mapsto Y^x_t$ has finitely many jumps on any compact interval, and is constant between its jumps, and if $0<s<t$,
 \begin{align}\label{monot}
 \{x,\, Y^{x}_t\not=Y^{x-}_t\}\subset\{x,\, Y^{x}_s\not=Y^{x-}_s\}.
 \end{align}
Let us now  fix $\delta, M>0$. For almost any realization of $Y$, the mapping $x \mapsto Y_\delta^x$ has only finitely many jumps on $(0,M]$. Let $n$ be so large that there is at most one of those jumps in each interval $(k2^{-n}, (k+1)2^{-n}]$, for $k\le M2^n-1$. 
Then for each $x$ that belongs to an interval $(k2^{-n}, (k+1)2^{-n}]$ which contains no jump of $x \mapsto Y_\delta^x$, and
 for any $n'>n$, we have $Z^{n',x}_t=Z^{n,x}_t$ for any $t\ge\delta$.
 
 Since $\delta$ and $M$ are arbitrary positive reals, we have shown that
 \begin{align}\label{newZ}
 \tilde Z^x_t:=\text{ a.s. }\lim_{n\to\infty}Z^{n,x}_t
 \end{align} exists for all $t\ge0$, $x>0$. The thus constructed random field $\{\tilde Z^x_t,\, t\ge0, x>0\}$ has the same law as the solution of the SDE \eqref{eqDL}, and satisfies \eqref{comparStrong} and hence also  
  \begin{align}\label{jumpsZY}
 \{x,\, \tilde Z^{x}_t\not=\tilde Z^{x-}_t\}\subset\{x,\, Y^{x}_t\not=Y^{x-}_t\}
 \end{align} 
 for all $t>0$. 
We still have to show that $\tilde Z$ satisfies \eqref{eq:DLcoupled}. It is plain that for any $\delta>0$,
 $$\tilde Z^x_t=\tilde Z^x_\delta+\int_\delta^t f(\tilde Z^x_s)ds+2\int_\delta^t\int_{A^x_s(\tilde Z)}W(ds,d\xi).$$
 In order to deduce that $\tilde Z$ satisfies \eqref{eq:DLcoupled}, it remains to show that $\tilde Z^x_\delta\to x$ a.s., as $\delta\to0$,
 which follows readily from the equality of the laws of $\tilde Z$ and $Z$.
 \epf

 \begin{corollary}\label{corRC}
 For any $t>0$, $x\mapsto Z^x_t$ has finitely many jumps 
 on any compact interval, and is constant between these jumps. 
  \end{corollary}
 \bpf 
The assertion follows from the fact that $\tilde Z$ possesses that property, as a consequence of \eqref{comparStrong} and the properties of $Y$.
\epf

  From the properties of the map $x\mapsto Z^x$, we infer that $x\mapsto \zeta^x:=\zeta(Z^x)$ is increasing and right continuous, constant between its jumps,
 with a.s. finitely many jumps on any compact subinterval of $(0,+\infty)$, and a sufficient condition  is given in \cite{LP} for the limit $\zeta^\infty$ to be a.s. finite.
 
 We have moreover
\begin{corollary}\label{cor2}
For any $s>0$,
 $$\P\left(\bigcup_{t>s}\{x,\ Z^x_t\not=Z^{x-}_t\}\subset\{x,\ Z^x_s\not=Z^{x-}_s\}\ \text{for all }x>0\right)=1.$$
\end{corollary}
\bpf
 Let us first fix $t>s$ and $x>0$. We have
 $$Z^x_t-Z^{x-}_t=Z^x_s-Z^{x-}_s+\int_s^t[f(Z^x_r)-f(Z^{x-}_r)]dr+2\int_s^t\int_{Z^{x-}_r}^{Z^x_r}W(dr,d\xi).$$
 Consequently, taking the conditional expectation given $Z^x_s-Z^{x-}_s$, and using both \eqref{1} and Gronwall's Lemma, we obtain
 $$\E[Z^x_t-Z^{x-}_t| Z^x_s-Z^{x-}_s]\le [Z^x_s-Z^{x-}_s]\exp(\theta (t-s))\quad \mbox{a.s.}$$
 This shows that $Z^x_t=Z^{x-}_t$ a.s. on the event $\{Z^x_s=Z^{x-}_s\}$. 
 
For all $M>0$ it follows from \eqref{jumpsZY} and Theorem \ref{def-coupling} that  $\{0<x\le M;\ Z^x_s\not=Z^{x-}_s\}$ is a.s. a random finite set. Let $0< V_1<V_2$ be $\sigma\{Z^y_r,\ y>0, r\le s\}$-measurable and such that $Z^x_s-Z^{x-}_s=0$ for all $V_1\le x<V_2$. From the above argument, for any $x\in[V_1,V_2)$, $t\ge s$,
 $Z^x_t-Z^{x-}_t=0$ a.s. Since $Z^x_t-Z^{x-}_t$ is continuous in $t$ and right continuous in $x$, 
 $$\P(Z^x_t-Z^{x-}_t=0,\ \text{for all }V_1\le x<V_2,\ t\ge s)=1.$$
 The result follows from the fact that the set $\{0<x\le M;\ Z^x_s-Z^{x-}_s=0\}$ is a.s. a finite union of intervals of the form $[V_1,V_2)$,
 and $M>0$ was arbitrary.
 \epf

  \begin{remark}\label{notuse}
  We believe that the coupling  constructed in Theorem \ref{def-coupling} is interesting in its own right.
 In the rest of this paper
  we shall exploit  its Corollary  \ref{corRC}.
 \end{remark}
 
 \section{An SDE for the path-valued Markov process}\label{sec:3}
Again, let  $Z^x$ be the solution of \eqref{eqDL}, with $f$ satisfying conditions \eqref{1}, \eqref{2} and \eqref{3}.
 From now on, the process  $Y^x$ will be the solution of 
\begin{equation}\label{eq:Y0}
Y^x_t=x+2\int_0^t\int_0^{Y^x_s}W(ds,d\xi).
\end{equation} 
 We shall use the notation
 $$F(a,b)=f(a+b)-f(a).$$
 
 Let $x,y>0$, and define 
 \begin{equation}\label{UV}
 V^{x,y}_t=Z^{x+y}_t-Z^x_t,  \quad U^{x,y}_t=Y^{x+y}_t-Y^x_t, \quad
 t\ge0.
 \end{equation} 
 We can couple these stochastic processes, by representing them as solutions
 \begin{align}
 V^{x,y}_t&=y +\int_0^t\ F(Z^x_s,V^{x,y}_s)ds
 +2\int_0^t\int_{Y^x_s}^{Y^x_s+V^{x,y}_s} W(ds,d\xi),\label{eqV} \\
 U^{x,y}_t&=y
 +2\int_0^t\int_{Y^x_s}^{Y^x_s+U^{x,y}_s} W(ds,d\xi), \label{eqU}
 \end{align}
with a $W$ different from (but having the same distribution as) the one appearing in \eqref{eqDL} and   \eqref{eq:Y0}, and  leading to a pair $(U,V)$ that has the same marginal distributions as the ones specified by \eqref{UV}.

We now define a Girsanov--Radon--Nikodym derivative, which will play an essential role in the sequel. For $z\in E$, $t>0$ and $U$ as in \eqref{eqU} (or, as we will need it later, also for some other $\R_+$--valued continuous semimartingale $U$ with quadratic variation $d \langle U\rangle_s = 4 U_s ds$), we define 
\begin{equation}\label{GRND}
L_t(z,U)=\exp\left(\frac{1}{4}\int_0^t\frac{F(z(s),U_s)}{U_s}dU_s
-\frac{1}{8}\int_0^t\frac{F(z(s),U_s)^2}{U_s}ds\right),
\end{equation}
where we use the convention $\frac{F(z,0)}{0}=0$. It follows from \eqref{2} that
for any $M <\infty$, $a\in [0,M]$ and $b \in [0,1]$,
\begin{align}\label{boundF}
 \frac{|F(a,b)|}{\sqrt b} \,  \le C_M;
\end{align} 
hence $L_t(z,U)$ is a well--defined random variable. 
We shall also consider  $L_t(Z,U)$, where $z$ is replaced by the process $Z$, solution of \eqref{eqDL} with some initial condition $Z_0=x$. Note that whenever we consider 
$L_t(Z,U)$, the processes $Z$ and $U$ will always be  mutually independent.

Finally $L(Z,U)$ (resp. $L(z,U)$)
will be defined by
\begin{equation}\label{LZU}
L(Z,U)=L_\infty(Z,U)=L_\zeta(Z,U)\quad  (\mbox{resp. }L(z,U)=L_\infty(z,U)=L_\zeta(z,U)),
\end{equation}
where $\zeta=\zeta(U) = \inf\{t>0,\ U_t=0\}$.
We shall consider the r.v. $L(Z,U)$ (or $L(z,U)$) only when $\zeta<\infty$ a.s., which e.g. is
the case if $U$ solves  \eqref{eqU}; hence the above quantities are well defined.

 We have 
 \begin{proposition}\label{GirsVU}
 For $V$ and $U$ as in \eqref{eqV}, \eqref{eqU}, the law of $\{V^{x,y}_t,\ 0\le t\le \zeta\}$ is absolutely continuous with respect to the law of
 $\{U^{x,y}_t,\ 0\le t\le \zeta\}$, and the Radon--Nikodym derivative is 
 $L(Z^x,U^{x,y})$. 
 \end{proposition}
 \bpf
 For simplicity, we suppress the superindices $x$ and $y$. We consider the filtration $\{\mathcal F_t, t \ge 0\}$ defined by $\mathcal F_t := \sigma ((Z_s,U_s) : 0\le s\le t\wedge \zeta)$ and introduce the local martingale
 $$L_t=\exp\left(\frac{1}{4}\int_0^{t\wedge\zeta} \frac{F(Z_s,U_s)}{U_s}dU_s-\frac{1}{8}\int_0^{t\wedge\zeta}\frac{F(Z_s,U_s)^2}{U_s}ds\right),$$
 where again $\zeta=\inf\{t>0,\ U_t=0\}$ is  the extinction time of $U$. Define, for each $n\ge1$,
 $T_n=\inf\{t>0,\ \int_0^tU_s^{-1}F^2(Z_s,U_s)ds>n\}\wedge\zeta$. It is plain that the sequence of events 
 $A_n=\{T_n=\zeta\}$ is increasing. 
Thus from the fact that $\zeta<\infty$ $\P$ a.s. together with  the assumption \eqref{2} it follows that $\P(\bigcup_n A_n)=1$. Moreover 
 for any fixed $n\ge1$, $(L_{t\wedge T_n})_{t\ge0}$ is a uniformly integrable martingale,
 and if we define $\Q_n$ on $\F_{T_n}$ by
 $$\frac{d\Q_n|_{\F_{T_n}}}{d\P|_{\F_{T_n}}}=L_{T_n},$$
 we have that the law of $(U_{t\wedge T_n})_{t\ge0}$ under $\Q_n$ equals the law of the process $(V_{t\wedge T_n})_{t\ge 0}$.
 It follows (see e.g. Proposition 3.5 in \cite{BP}) that there exists a unique probability measure $\Q$ on $\F_\zeta=\sigma(\cup_n\F_{T_n})$ such that, 
for each $n\ge1$, its restriction to $\F_{T_n}$ coincides with $\Q_n$.
 It remains to show that $\Q \ll \P$, and that
 $$\frac{d\Q}{d\P}=L_{\zeta}.$$
 For this purpose, let $A\in\F_\zeta$ and $n\ge1$. Clearly $A\cap A_n\in\F_{T_n}$ and
 \begin{align*}
 \Q(A\cap A_n)&=\E_\P({\bf1}_{A\cap A_n}L_{T_n})\\
 &=\E_\P({\bf1}_{A\cap A_n}L_\zeta).
 \end{align*}
We have not only $\P(\bigcup_n A_n)=1$, but also $\Q(\bigcup_nA_n)=1$ (indeed condition \eqref{3} implies that $Z^{x+y} = Z^x+V^{x,y}$ goes extinct in finite time a.s., hence also $V=V^{x,y}$ has this property). Thus, by letting $n\to\infty$  in the above equality, we deduce from the monotone convergence theorem that
 $$\Q(A)=\E_\P({\bf1}_{A}L_{\zeta}).$$
 Since $L_\zeta= L(Z^x,U^{x,y})$, the proposition is proved.  \epf
 
 Let us define for each $x>0$ the sigma--field $\G^x=\sigma\{Z^\xi_t,\, 0<\xi\le x, \, t\ge0\}$. As a corollary of Proposition \ref{GirsVU} and of the independence of $Z^x$ and $U^x$ we obtain
 for an arbitrary $t>0$,
 \begin{align}\label{GirsUV}
 \E(V_t^{x,y}|\G^x)=\E(L(Z^x,U^{x,y})U^{x,y}_t|\G^x).
 \end{align}

 In order to achieve the goal of deriving an SDE for the path-valued process ${Z^x}, x \ge 0$, and in view of Proposition \ref{GirsVU}, we want to take the limit as $y\to0$ in the expression
 $$\frac{1}{y}\E(L(Z^x,U^{x,y})U^{x,y}_t|\G^x).$$
Note that the law of $U^{x,y}$ is  that of the unique
solution of the SDE
\begin{equation}\label{CFD}
U_t=y+2\int_0^t\sqrt{U_s}\, dB_s.
\end{equation}
 In particular that law (which is a probability measure on $E$) does not depend on $x$; we denote it by $\mathbf P_y$.
For any $t>0$, $A\in\sigma\{U_r,\ r\ge0\}$, $y>0$, 
 let  \begin{equation}\label{defQ}
 \mathbf Q_{y,t}(A)=y^{-1}\mathbf E_y(U_t;A).
 \end{equation}
 For $t > 0$ we now write $\F_t := \sigma\{U_r,\ 0\le r\le t\}$. The next result is known, tracing back to \cite{LN} Theorem 1 (see also \cite{RR}, Theorem 2, \cite{zL} Theorem 4.1 , and \cite{aL} Theorem 4.1 for more general versions). We give a proof, which is short, for the convenience of the reader.
 \begin{proposition}\label{fellim}
 For any fixed $t > 0$ and $y>0$,  the process $U =\{U_r, 0\le r\le t \}$ is, under $\mathbf Q_{y,t}$, a Feller process with immigration. More precisely, $U$ solves under $\mathbf Q_{y,t}$ the SDE
\begin{equation}\label{Fimm}
U_r=y+4r+2\int_0^r\sqrt{U_s}d{\bar B}_s,\ 0\le r\le t,
\end{equation}
where $\bar B$ is a $\mathbf Q_{y,t}$--standard Brownian motion. 
 \end{proposition}
 \bpf Denoting by $\zeta$  the extinction time of $U$, we have
 \begin{align*}
 \frac{U_{t\wedge\zeta}}{y}&=\exp(\log U_{t\wedge\zeta}-\log U_0),\\
 \log U_{t\wedge\zeta}-\log U_0&=2\int_0^{t\wedge\zeta}\frac{dB_s}{\sqrt{U_s}}-2\int_0^{t\wedge\zeta}\frac{ds}{U_s},
 \end{align*}
 hence
 $$\frac{d\mathbf Q_{y,t}|_{\F_t}}{d\mathbf P_y|_{\F_t}}=\exp\left( 2\int_0^{t\wedge\zeta}\frac{dB_s}{\sqrt{U_s}}-2\int_0^{t\wedge\zeta}\frac{ds}{U_s}  \right),$$
 which by Girsanov's theorem implies the result, if we let $\bar B_r=B_r-2\int_0^r(U_s)^{-1/2}ds$.
 
 Note that we can apply Girsanov's theorem here, since $\E_{\mathbf P_y} U_t=y$ implies that
 $$\E_{\mathbf P_y}\left[\exp\left(2\int_0^{t\wedge\zeta}\frac{dB_s}{\sqrt{U_s}}-2\int_0^{t\wedge\zeta}\frac{ds}{U_s}\right)\right]=1.$$\epf
 \\
 As a corollary to Proposition~\ref{fellim} (and the Markov property) we get that  for $t > 0$, the process $\{U_r, r\ge 0 \}$ under $\mathbf Q_{y,t}$ solves the SDE
\begin{equation}\label{FDimmt}
U_r=y+4\; (r\wedge t)+2\int_0^r\sqrt{U_s}d\bar B_s,\quad r\ge0.
\end{equation}
  It is immediate  that the limit  of $\mathbf Q_{y,t}$ exists as $y\to 0$. We will denote this limit by $\mathbf Q_{0,t}$, and note that
it is the law of
 $\{U_r,\ r \ge 0\}$,  the solution of 
 $$U_r=4(r\wedge t)+2\int_0^r\sqrt{U_s}d\bar B_s,\ r \ge 0.$$
 For $y\ge0$, we denote by $\mathbf Q_{y,\infty}$ the law of the process $U = \{U_r, r \ge 0\}$, that satisfies \eqref{FDimmt} with $t=\infty$. It is well known (see e.g. \cite{aL})  that $\mathbf Q_{y,\infty}$ is the law of a critical Feller process conditioned to never die out, hence
 for $y >0$ we have
\begin{align}\label{Upos}
\mathbf Q_{y,\infty}(U_s>0\quad  \forall s \ge 0)=1,
\end{align}
whereas for $y=0$ we observe that
\begin{align}\label{Uspos}
\mathbf Q_{0,\infty}(U_s>0 \quad \forall s > 0 )=1.
\end{align}
 For $t< \infty$, $r > 0$ and $z\in E$,  from $L_r(z,U)$ and $L(z,U)$ defined as in \eqref{GRND} and \eqref{LZU} we obtain the $\mathbf Q_{y,t}$-a.e. defined measurable functions  $u \mapsto L_r(z,u)$ and $u \mapsto L(z,u)$. Under $\mathbf{Q}_{y,\infty}$, 
 we shall consider only $L_r(z,U)$, since $\zeta=+\infty$ \mbox{$\mathbf{Q}_{y,\infty}$ a.s.}

 For any $t>0$ and any partition $\mathcal P = \{x_k, \, k = 0,1,\ldots\}$ with $0=x_0<x_1<\cdots$ we have that
 \begin{align*}
 Z^{x_\ell}_t&=\sum_{k=1}^\ell(x_k-x_{k-1})\E\left(\frac{Z^{x_k}_t-Z^{x_{k-1}}_t}{x_k-x_{k-1}}\Big|Z^{x_{k-1}}\right)+M^{x_\ell, \mathcal P}_t\\
 &=\sum_{k=1}^\ell(x_k-x_{k-1})\int_EL(Z^{x_{k-1}},u)\mathbf Q_{x_k-x_{k-1},t}(du)+M^{x_\ell, \mathcal P}_t,
 \end{align*}
 where we have used \eqref{GirsUV} and \eqref{defQ}, and where we define
 $$M^{x_\ell, \mathcal P}_t=\sum_{k=1}^\ell\left[ Z^{x_k}_t-Z^{x_{k-1}}_t-\E(Z^{x_k}_t-Z^{x_{k-1}}_t|\G^{x_{k-1}}) \right].$$
 
 For every $n \in \mathbb N$, we consider the partition  $\mathcal P_n:= \{x^n_k=k2^{-n}, \, k\ge1\}$. It follows from the above arguments that if $x$ is a dyadic number and
 $n$ is large enough, then, with $M^{x,n}:=M^{x,\mathcal P_n}$,
  \begin{equation}\label{eq:tn}
  Z^x_t=\sum_{k=1}^{x2^n}2^{-n}\int_EL(Z^{(k-1)2^{-n}},u)\mathbf Q_{2^{-n},t}(du)+M^{x,n}_t.
   \end{equation}
Our aim is to show   convergence of the right hand side as  $n\to\infty$, leading to
 \begin{equation}\label{eq:t}
 Z^x_t=\int_{[0,x]\times E}L(Z^\xi,u)\mathbf Q_{0,t}(du)d\xi+M^x_t,
 \end{equation}
 where $\{M^x_t,\ x\ge0\}$ is a $\G^x$--martingale.
 To this purpose we start by proving
 \begin{lemma}\label{le:Q}
 For any $y\ge0$, $t>0$  and $z \in E$,
  \begin{equation}\label{Qlemma}
  \int_E L(z,u)\mathbf Q_{y,t}(du)= \int_E L_t(z,u)\mathbf Q_{y,\infty}(du).
  \end{equation}
  \end{lemma}
  \bpf 
Since $\mathbf Q_{y,\infty}$ and $\mathbf Q_{y,t}$ coincide when restricted to $\mathcal F_t$, and since $L_t(z,\cdot)$ is $\mathcal F_t$-measurable, the right  hand side of \eqref{Qlemma}  equals 
  \begin{equation}\label{intLQ}
  \int_E L_t(z,u)\mathbf Q_{y,t}(du).
  \end{equation}
   It thus remains to show that this is equal to the left hand side of \eqref{Qlemma}. Because of \eqref{FDimmt}, 
 the process $(U_s)_{s\ge t}$  under $\mathbf Q_{y,t}$ is a driftless Feller diffusion. Thus the exponential martingale 
 $$(M_s)_{s\ge t}:= (L_s(z,U)/L_t(z,U))_{s\ge t}$$
  constitutes a process of Girsanov densities with respect to $\mathbf Q_{y,t}$, with the property that $U$ under the transformed measure satisfies the SDE  
  \begin{equation}\label{eq:V}
  dU_s = F(z_s,U_s)ds +2\sqrt{U_s} d\tilde B_s,\quad s \ge t.
  \end{equation}
  for a standard Brownian motion $\tilde B$. From this we conclude in precisely the same manner as in the proof of Proposition \ref{GirsVU} that $\mathbb E_{\mathbf Q_{y,t}} [M_\zeta] = 1$, where $\zeta$ denotes the extinction time of $U$. This implies  that $\mathbb E_{\mathbf Q_{y,t}} [M_\zeta|\F_t] = 1$ a.s., hence 
  $\mathbb E_{\mathbf Q_{y,t}}[L(z,U)]= \mathbb E_{\mathbf Q_{y,t}} [L_\zeta(z,U)] = E_{\mathbf Q_{y,t}} [L_t(z,U)]$, that is, the l.h.s. of \eqref{Qlemma} indeed equals \eqref{intLQ}.  \epf
 
 For $y \ge 0$, $z\in E$ and $\tilde B$ a standard Brownian motion, consider the SDE
 \begin{equation}\label{FLimm}
V_r=y+4r+\int_0^rF(z_s,V_s)ds + 2\int_0^r\sqrt{V_s}d\tilde B_s,\ r \ge 0,
\end{equation}
 and denote the law of its solution by $\tilde {\mathbf Q}^z_{y,\infty}$. 
      \begin{lemma}\label{le:expression}
 For any $y\ge0$, $t>0$ and $z \in E$,
 \begin{equation*}
 \int_E L_t(z,u)\mathbf Q_{y,\infty}(du)= \int_E \exp\left({\int_0^tu_s^{-1}F(z_s,u_s)\,ds}\right)\tilde {\mathbf Q}^z_{y,\infty}(du).
\end{equation*}
  \end{lemma}
 \bpf  Let $y\ge0$, $t>0$ and $z \in E$ be fixed throughout the proof. We recall that $\mathbf Q_{y,\infty}$ has been defined as the law of the solution $U$ of the SDE
\begin{equation}\label{FDI}
 U_r=y+4r +2\int_0^r\sqrt{U_s}d\bar B_s,\, r\ge 0.
 \end{equation}
We note that
\begin{equation}\label{Lrep}
L_t(z,U)=G_t(z,U)\ e^{\int_0^tU_s^{-1}F(z_s,U_s)ds},
\end{equation}
with $G_t$ being defined as 
$$G_t = G_t(z,U) = \exp\left(\frac{1}{2} \int_0^t\frac{F(z_s,U_s)}{\sqrt {U_s}} d\bar B_s - \frac{1}{8} \int_0^t\frac{F^2(z_s,U_s)}{U_s}\, ds\right).$$
It is plain that $(G_r, r\ge 0)$ is a local martingale under $\mathbf Q_{y,\infty}$. In order to conclude that ${\tilde {\mathbf Q}}^z_{y,\infty}|\mathcal F_t$ has density $G_t$ w.r.to $\mathbf Q_{y,\infty}|\mathcal F_t$, and thus to infer the assertion of the lemma,  it will be sufficient to ensure the Girsanov condition
\begin{equation}\label{GCt}
\mathbb E_{\mathbf Q_{y,\infty}}[G_t] = 1 . \end{equation}
For this we proceed similarly as in the proof of Proposition \ref{GirsVU} and define for $n=0,1,2, \ldots$
$$T_n:= \inf\{r: \int_0^rU_s^{-1}F^2(z_s,U_s)\, ds \ge n\}\; ,\quad T_\infty:=\lim_{n\to\infty}T_n.$$
Since $(G_{t\wedge T_n}, n=0,1,2,\ldots)$ is a martingale, we can define a measure
$\tilde {\mathbf Q}_{y}$ on $\mathcal F_{t\wedge T_\infty}$ whose restriction to  $\mathcal F_{t\wedge T_n}$ is given by
$$\frac {d\tilde {\mathbf Q}_{y} | \mathcal F_{t\wedge T_n}}{d \mathbf Q_{y,\infty} | \mathcal F_{t\wedge T_n}} = G_{t \wedge T_n}.$$
Then  for all  $n \ge 1$, since $\{T_n\ge t\} \in \mathcal F_{t\wedge T_n}$, we have
$$\tilde {\mathbf Q}_{y}(\{T_n\ge t\}) = \mathbb E_{\mathbf Q_{y,\infty}}(\mathbf 1_{\{T_n\ge t\}}G_{t\wedge T_n}) =  \mathbb E_{\mathbf Q_{y,\infty}}(\mathbf 1_{\{T_n\ge t\}}G_{t}).$$
Letting $n\to \infty$ we obtain
\begin{equation}\label{GL}
\tilde {\mathbf Q}_{y}(\{T_\infty\ge t\})  =  \mathbb E_{\mathbf Q_{y,\infty}}(\mathbf 1_{\{T_\infty\ge t\}}G_{t}).
\end{equation}
Because of 
$$\{T_n\ge t\} = \{T_n\wedge t=t\} = \{\int_0^tU_s^{-1}F^2(z_s,U_s) ds\le n\},$$
we have
\begin{equation}\label{twoevents}
\{T_\infty \ge t\} =  \{\int_0^tU_s^{-1}F^2(z_s,U_s) ds< \infty\}.
\end{equation}
Under the measure $\mathbf Q_{y,\infty}$ the process $U$ solves  \eqref{Fimm}, hence it does not explode in finite time. Thus, from \eqref{2},
$$\mathbf Q_{y,\infty}(T_\infty \ge t) = 1.$$
On the other hand,  under the measure $\tilde {\mathbf Q}_{y}$,  the process  $U$ is the solution of
\begin{equation}\label{SDV}
 U_r=y+4r+\int_0^rF(z_s,U_s)ds+2\int_0^r\sqrt{U_s}d\tilde B_s
 \end{equation}
 up to the time   $T_\infty \wedge  t$.
Thanks to \eqref{1} the solution of \eqref{SDV} does not explode in finite time, consequently we deduce again from \eqref{2} that
$$\tilde {\mathbf Q}_{y} (T_\infty \ge t) = 1.$$
Thus \eqref{GL} simplifies to  \eqref{GCt}. 
 \epf
 
 We now prove the
 \begin{proposition}\label{prop:cont_in_t}
  For all $z\in E$, $y\ge0$, the mapping $t\mapsto \varphi(t):=\int_E L_t(z,u)\mathbf Q_{y,\infty}(du)$ is continuous on $(0,\infty)$,  
 satisfies $0\le\varphi(t)\le \exp(\theta t)$, and it converges to $1$ as $t \to 0$.
 \end{proposition} 
\bpf 
 Fix $z\in E$ and  $y\ge0$. Let $U$ be an $E$-valued random variable with distribution $\mathbf Q_{y,\infty}$. Then, for all $z \in E$, the mapping $t\to L_t(z,U)$ is a.s. continuous.
 From \eqref{Lrep} and \eqref{1},  we infer that $0\le L_t(z,U)\le G_t(z,U) e^{\theta t}$, hence for any $b>0$, the uniform integrability of the family of random variables
 $\{L_t(z,U),\ 0\le t\le b\}$  follows from that of
 $\{G_t (z,U),\ 0\le t\le b\}$, which in turn follows from its martingale property established in Lemma \ref{le:expression}. 
 This implies the claimed continuity. The fact that the integral equals 1 at $t=0$ follows from the fact that $L_0(z,u)=1$.
Finally, the fact that $\int_E L_t(z,u)\mathbf Q_{y,\infty}(du)\le\exp(\theta t)$ follows readily from Lemma \ref{le:expression}
 and \eqref{1}. 
 \epf
\begin{proposition}\label{prop:cont_in_z}
For any $t, y\ge0$, the mapping
$$z\mapsto\int_EL_t(z,u)\mathbf Q_{y,\infty}(du)$$
is continuous from $E$ into $[0,e^{\theta t}]$.
\end{proposition}
\bpf
(a) Let us first check that the mappings  $z\mapsto G_t(z,.)$ and $z\mapsto L_t(z,.)$ and  both are continuous in $\mathbf Q_{y,\infty}$-probability.  Clearly, for any $s >0$ the mapping
$z\mapsto(U_s^{-1}F(z(s),U_s),U_s^{-1/2}F(z(s),U_s))$ is $\mathbf Q_{y,\infty}$ a.s. continuous from
$E$ into $\R^2$. 
\\
(i) We start by considering the case $y > 0$, and recall \eqref{Upos} as well as \eqref{boundF}.
 Hence, if  $z_n(s)\to z(s)$ uniformly in $s\in [0,t]$, then 
\begin{align}\label{as1}
\int_0^t\left|U_s^{-1/2}F(z(s),U_s)-U_s^{-1/2}F(z_n(s),U_s)\right|^2ds\to 0 \quad \mathbf Q_{y,\infty} \mbox{ a.s. }
\end{align}
and
\begin{align}\label{as2}
\int_0^t U_s^{-1}F(z_n(s),U_s)ds\to \int_0^t U_s^{-1}F(z(s),U_s)ds \quad \mathbf Q_{y,\infty} \mbox{ a.s. }
\end{align}
(Note that \eqref{as2} follows from \eqref{boundF} and dominated convergence, since due to \eqref{Upos}  we have $\int_0^t U_s^{-1/2} ds < \infty\quad \mathbf Q_{y,\infty}$ a.s.)
From \eqref{as1} and \eqref{as2} we conclude that
$G_t(z_n,U)\to G_t(z,U)$ and  $L_t(z_n,U)\to L_t(z,U)$, both in $\mathbf Q_{y,\infty}$-probability.
\\
(ii)  It remains to treat the case $y=0$. Then $U_0=0$, and \eqref{Uspos} holds.
 Let $T=\inf\{0\le s\le 1,\ U_s\ge1\}$. Consider again a sequence $z_n$ in $E$ such that $z_n(s)\to z(s)$ uniformly in $[0,t]$.
We let $M=\sup_{n\ge1}\sup_{0\le s\le t}z_n(s)$, and $C_{M}$ be the associated constant appearing in 
\eqref{2}. Then whenever $0 < s\le T\wedge t$, 
$$|U_s^{-1/2}F(z_n(s),U_s){\bf1}_{\{U_s\le 1\}}|\le C_{M},$$
hence by bounded convergence
\begin{align}\label{as3}
\int_0^{t\wedge T}\left|U_s^{-1/2}F(z(s),U_s)-U_s^{-1/2}F(z_n(s),U_s)\right|^2ds\to 0 \quad \mathbf Q_{0,\infty} \mbox{ a.s. }
\end{align} as $n\to\infty$. The convergence \eqref{as1} for $y=0$ now follows from \eqref{as3} and the above part 1 of the proof, with the strong Markov property applied to the stopping time $T$ and $y=1$. 
  From this we conclude that
$G_t(z_n,U)\to G_t(z,U)$ in $\mathbf Q_{0,\infty}$-probability.
To obtain \eqref{as2} from \eqref{boundF} and dominated convergence also in the case $y=0$,  we observe, using \eqref{Fimm}
 and It\^o's formula applied to $\sqrt{U_s}$,
  that $$\int_0^{t} U_s^{-1/2} ds =  \frac{2}{3}[\sqrt{U_t} -B_{t} ]< \infty \quad \mathbf Q_{0,\infty} \text{ a.s.}$$
On the other hand from \eqref{boundF} and \eqref{Uspos} we have for all $s>0$
$$\frac{| F(z_n(s), U_s)|}{U_s}{\bf1}_{\{U_s\le1\}}\le \frac{C_M}{\sqrt{U_s}} \quad \mathbf Q_{0,\infty} \text{ a.s.},$$  thus we conclude by Lebesgue's dominated convergence theorem that
\begin{align*}
\int_0^{t\wedge T} U_s^{-1}F(z_n(s),U_s)ds\to \int_0^{t\wedge T} U_s^{-1}F(z(s),U_s)ds \quad \mathbf Q_{0,\infty} \mbox{ a.s. }
\end{align*}
This together with \eqref{as1} for $y=1$  implies \eqref{as1} also for $y=0$. 

We have thus established  \eqref{as1} and  \eqref{as2} for $y \ge 0$, which yields the claimed continuity of the mapping $z\mapsto L_t(z,.)$ in $\mathbf Q_{y,\infty}$-probability. 
\\
(b) In order to conclude the proof we need  uniform integrability of the family $L_t(z_n,.), \, z \in E$, with respect to $\mathbf Q_{y,\infty}$, where $z_n, z \in E$ and $z_n \to z$.  This can be established in the very same manner as we did in the proof of Proposition \ref{prop:cont_in_t}, after observing that $0\le L_t(z_n,.)\le G_t(z_n,.) e^{\theta t}$, and that $G_t(z_n,.)$  not only converges in  $\mathbf Q_{y,\infty}$-probability towards $G_t(z,.)$ due to part (a), but also is uniformly integrable because of $\mathbb E_{\mathbf Q_{y,\infty}} [ G_t(z_n,U)] = 1$ for all $n\ge1$, and $\mathbb E_{\mathbf Q_{y,\infty}} [ G_t(z,U)] = 1$.
\epf

Combining this result with Lemma \ref{cadlag}, we deduce
\begin{corollary}\label{corintRC}
For any $t, y\ge0$, the mapping
$$x\mapsto\int_EL_t(Z^x,u)\mathbf Q_{y,\infty}(du)$$
is a.s. c\`adl\`ag from $\R_+$ into $\R_+$.
\end{corollary}
We are now in a position to establish
\begin{proposition}\label{prop:conv}
For any $t,x\ge0$, 
$$\sum_{k=1}^{x2^n}2^{-n}\int_EL_t(Z^{(k-1)2^{-n}},u)\mathbf Q_{2^{-n},\infty}(du)\to 
\int_{[0,x]\times E}L_t(Z^\xi,u)\mathbf Q_{0,\infty}(du)d\xi$$
in probability, as $n\to\infty$.
\end{proposition}
\bpf  
We first show that
\begin{align}\label{firstconv}
\sum_{k=1}^{x2^n}2^{-n}\int_EL_t(Z^{(k-1)2^{-n}},u)\left[\mathbf Q_{2^{-n},\infty}(du)-\mathbf Q_{0,\infty}(du)\right]\to0
\end{align}
as $n\to\infty$. To this purpose we define for each $n\ge1$ and $z\in E$
$$H_n(z)=\int_EL_t(z,u)\left[\mathbf Q_{2^{-n},\infty}(du)-\mathbf Q_{0,\infty}(du)\right].$$
It follows from Lemma \ref{le:expression} that 
\begin{equation}\label{eqstar}
\begin{split}
H_n(z)=
\E&\left[\exp\left(\int_0^t(\mathcal V^{2^{-n}}_s(z))^{-1}F(z_s,\mathcal V^{2^{-n}}_s(z))ds\right) \right. \\
&\quad \left.-\exp\left(\int_0^t(\mathcal V^0_s(z))^{-1}F(z_s,\mathcal V^0_s(z))ds\right)\right],
\end{split}
\end{equation}
where $\mathcal V^{2^{-n}}(z)$ (resp. $\mathcal V^0(z)$) denotes the solution of the SDE \eqref{FLimm} with $y=2^{-n}$ (resp. with $y=0$).
For $\xi \in [0,x]$ we put
$$h_n(\xi)=\sum_{k=1}^{x2^n}H_n(Z^{(k-1)2^{-n}}){\bf1}_{[(k-1)2^{-n}, k2^{-n})}(\xi).$$
Whenever $\xi\in[(k-1)2^{-n}, k2^{-n})$, we briefly write 
\begin{align}\label{xin}
\xi_n:=[\xi2^{n}]2^{-n}=(k-1)2^{-n},\end{align}
hence as $n\to\infty$,
$\xi_n\to\xi$ and $Z^{\xi_n}\to Z^{\xi-}$ a.s. Also, $h_n$ can be rewritten as 
$$h_n(\xi) = H_n(Z^{\xi_n}).$$
From \eqref{1}, the expression in the expectation on the right hand side of \eqref{eqstar}  is bounded in absolute value by $2\exp(\theta t)$. 
Hence we infer from 
Lemma \ref{le:conv} below and Lebesgue's dominated convergence theorem that for all $\xi > 0$
$$H_n(Z^{\xi_n}) \to 0  \mbox { in probability as } n\to \infty.$$
Since the left hand side of \eqref{firstconv} equals $\int_0^xh_n(\xi)d\xi$, and since  $|H_n(z)|\le 2\exp(\theta t)$, 
the assertion \eqref{firstconv} follows by dominated convergence.
It thus remains to show that
$$\sum_{k=1}^{x2^n}2^{-n}\int_EL_t(Z^{\xi_n},u)\mathbf Q_{0,\infty}(du)\to\int_{[0,x]\times E}
L_t(Z^\xi,u)\mathbf Q_{0,\infty}(du)d\xi$$
a.s, as $n\to\infty$. 
This follows readily from Corollary \ref{corintRC} together with the elementary fact that
for any
right--continuous mapping $\xi\mapsto\mathcal{A}(\xi)$  from $[0,+\infty)$ into $[0,1]$, and any dyadic $x>0$, as $n\to\infty$, one has
$$2^{-n}\sum_{k=1}^{[x 2^n]}\mathcal{A}((k-1)2^{-n})\to\int_{[0,x]}\mathcal{A}(\xi)d\xi.\qquad \qquad \mbox{\epf}$$
We finally establish the following result
\begin{lemma}\label{le:conv} Let $\xi > 0$ and $\xi_n$ be as in \eqref{xin}, $\mathcal V^{n,n}:=\mathcal V^{2^{-n}}(Z^{\xi_n})$,  $\mathcal V^n:=\mathcal V^0(Z^{\xi_n})$.
Then 
$$\int_0^t(\mathcal V^{n,n}_s)^{-1}F(Z^{\xi_n}_s,\mathcal V^{n,n}_s)ds-\int_0^t(\mathcal V^n_s)^{-1}F(Z^{\xi_n}_s,\mathcal V^n_s)ds\to0$$
in probability, as  $n\to\infty$.
\end{lemma}
\bpf
We first note that $\mathcal V^{n,n}_s>0$ a.s. for all $s\ge0$, while $\mathcal V^n_s>0$ for all $s>0$, but $\mathcal V^n_0=0$.
These facts follow from our assumption \eqref{2}, which implies that when those solutions get close to zero, their drift is bigger than 2.
Consequently, since for any $s$  the mapping $v\mapsto v^{-1}F(Z^{\xi_n}_s,v)$   is locally bounded and continuous away from $v=0$, we conclude that
for any $0<\delta\le1$, both
\begin{align*}
\int_\delta^t(\mathcal V^{n,n}_s)^{-1}F(Z^{\xi_n}_s,\mathcal V^{n,n}_s)ds&\to\int_\delta^t(\tilde{\mathcal V}_s)^{-1}F(Z^{\xi-}_s,\tilde{\mathcal V}_s)ds,\ \text{and}\\
 \int_\delta^t(\mathcal V^n_s)^{-1}F(Z^{\xi_n}_s,\mathcal V^n_s)ds&\to\int_\delta^t(\tilde{\mathcal V}_s)^{-1}F(Z^{\xi-}_s,\tilde{\mathcal V}_s)ds
\end{align*}
a.s. , where $\tilde {\mathcal V}$ denotes the solution of the SDE \eqref{FLimm} with $y=0$ and $z=Z^{\xi-}$.
Since for any $\eps>0$,
\begin{align*}
\P&\left(\left|\int_0^t(\mathcal V^{n,n}_s)^{-1}F(Z^{\xi_n}_s,\mathcal V^{n,n}_s)ds- \int_0^t(\mathcal V^n_s)^{-1}F(Z^{\xi_n}_s,\mathcal V^n_s)ds\right|>3\eps\right)\\
&\quad\le\P\left(\left|\int_\delta^t(\mathcal V^{n,n}_s)^{-1}F(Z^{\xi_n}_s,\mathcal V^{n,n}_s)ds- \int_\delta^t(\mathcal V^n_s)^{-1}F(Z^{\xi_n}_s,\mathcal V^n_s)ds\right|>\eps\right)
\\&\quad+\P\left(\left|\int_0^\delta (\mathcal V^{n,n}_s)^{-1}F(Z^{\xi_n}_s,\mathcal V^{n,n}_s)ds\right|>\eps\right)
\\&\quad+\P\left(\left|\int_0^\delta (\mathcal V^n_s)^{-1}F(Z^{\xi_n}_s,\mathcal V^n_s)ds\right|>\eps\right),
\end{align*}
the Lemma will follow from the fact that for any $\eps>0$, $\eta>0$, there exists $0<\delta\le1$ small
enough such that
\begin{equation}\label{eq:eta}
\begin{split}
\P\left(\left|\int_0^\delta (\mathcal V^{n,n}_s)^{-1}F(Z^{\xi_n}_s,\mathcal V^{n,n}_s)ds\right|>\eps\right)&\le \eta,\\
\P\left(\left|\int_0^\delta (\mathcal V^n_s)^{-1}F(Z^{\xi_n}_s,{\mathcal V_s^n})ds\right|>\eps\right)&\le \eta.
\end{split}
\end{equation}
It is plain that $\mathcal V^n_s\le \mathcal V^{n,n}_s\le \overline{V}_s$, where $\overline{V}$ solves the SDE
\begin{align}\overline{V}_t=1+(4+\theta)t+2\int_0^t\sqrt{\overline{V}_s}dB_s,\ t\ge0,
\end{align}
and for any $\eta>0$, there exists $M>0$ such that
$$\P\left(\sup_{0\le s\le \delta}\overline{V}_s>M\right)\le\frac{\eta}{2}.$$
In order to show \eqref{eq:eta} we consider the events  $\Omega_{M}=\{\sup_{0\le s\le 1,n\ge1}Z^{\xi_n}_s\le M\}$.
From \eqref{2} we infer $\sup_{0\le u\le M;\ 0\le z\le M}u^{-1/2}|F(z,u)|\le c_M$ for some finite constant $c_M$ depending on $M$. Thus, if we let $\tau_M=\inf\{s>0,\ \overline{V}_s> M\}$,
it follows from It\^o's formula that, writing $\mathcal V_s$ for either $\mathcal V^{n,n}_s$ or $\mathcal V^n_s$, 
$$\int_0^{\delta\wedge\tau_M}\frac{ds}{\sqrt{\mathcal V_s}}=\frac{2}{3}\left[\sqrt{\mathcal V_{\delta\wedge\tau_M}}-\sqrt{\mathcal V_0}-\frac{1}{2}\int_0^{\delta\wedge\tau_M}
\frac{F(Z^{\xi_n}_s,\mathcal V_s)}{\sqrt{\mathcal V_s}}ds-\tilde B_{\delta\wedge\tau_M}\right],$$
hence we have on $\Omega_{M}$
\begin{equation}\label{eq:estim}
0\le\int_0^{\delta\wedge\tau_M}\frac{ds}{\sqrt{\mathcal V_s}}\le \frac{2}{3}\left[\sqrt{\overline{ V}_{\delta\wedge\tau_M}}+\delta\frac{c_M}{2}
-\tilde B_{\delta\wedge\tau_M}\right].
\end{equation}
Now for both $\mathcal V=\mathcal V^{n,n}$ and $\mathcal V=\mathcal V^n$,
\begin{align*}
&\P\left(\left|\int_0^\delta \mathcal V^{-1}_sF(Z^{\xi_n}_s,\mathcal V_s)ds\right|>\eps\right)\\
&\le\P(\Omega_{M}) + \P(\tau_M<\delta) +
\P\left(\left|\int_0^{\delta\wedge\tau_M} \mathcal V_s^{-1}F(Z^{\xi_n}_s,\mathcal V_s)ds\right|>\eps; \Omega_{M}\right)\\
&\le \P(\Omega_{M}) + \frac{\eta}{2}+\P\left(\int_0^{\delta\wedge\tau_M} \frac{ds}{\sqrt{\mathcal V_s}}>c_M^{-1}\eps; \Omega_{M}\right).
\end{align*}
Finally, since $\mathbb P(\Omega_M) \to 0$ as $M \to \infty$, by choosing first $M$ sufficiently large and then $\delta$ sufficiently small, \eqref{eq:eta} follows readily from the previous estimate and \eqref{eq:estim}.
\epf

Guided by \eqref{eq:t} we now define
\begin{align}\label{Mx}
 M^x_t:=Z^x_t-\int_{[0,x]\times E}L_t(Z^\xi,u)\mathbf Q_{0,t}(du)d\xi, \quad t\ge 0, \, x\ge 0.
 \end{align}
 \begin{remark}\label{pathidentity}
 Due to Lemma \ref{le:Q} and Proposition \ref{prop:cont_in_t}, also the second summand on the r.h.s. of  \eqref{Mx} is a.s. continuous in $t$ for all $x$, hence \eqref{Mx} lifts to an identity for continuous-path valued processes (indexed by $x$). 
 \end{remark}
\begin{lemma}
 For any fixed $t>0$,  $\{M^x_t,\ x>0\}$, defined by \eqref{Mx},  is a c\`adl\`ag $(\G^x)$-martingale.
 \end{lemma}
 \bpf
 The first term on the right hand side of \eqref{Mx} is c\`adl\`ag, and the second is continuous. As for the martingale property,
 we first note that $M^x_t$ is integrable (since $Z^x_t$ is integrable, and  the last term in the right hand side of \eqref{Mx} takes values in $[0,xe^{\theta t}]$ due to Proposition \ref{prop:cont_in_t}) and $\G^x$-measurable for all $x>0$.
 We next show that if $0<a<x$  are two dyadic real numbers, then
 \begin{equation}\label{eq:mart}
 \E[M^x_t-M^a_t|\G^y]=0.
 \end{equation}
 Indeed, there exists $m\in\N$ such that $am$ and $xm$ are integers. Then for any $n\ge m$, 
 $M^{x,n}_t$ and $M^{a,n}_t$ defined by \eqref{eq:tn} satisfy
 \begin{equation}\label{eq:martn}
 \E[M^{x,n}_t-M^{a,n}_t|\G^y]=0.
 \end{equation}
  It follows from the above arguments that $M^{x,n}_t\to M^x_t$ and $M^{a,n}_t\to M^a_t$ a.s.
 as $n\to\infty$. Moreover $M^{x,n}_t$ is the difference of the integrable r.v. $Z^x_t$ which
 does not depend upon $n$, and a nonnegative r.v. which depends upon $n$ and is uniformly bounded by $xe^{\theta t}$. 
 Consequently the convergence holds also in $L^1$. The same is true for the sequence  $M^{a,n}_t$.
Hence \eqref{eq:mart} follows from \eqref{eq:martn}.

Suppose now that $x$ and $a$ are arbitrary positive real numbers, satisfying again $0<a<x$. Let $x_n$ (resp. $a_n$) be a decreasing sequence of dyadic reals, such that $x_n\to x$ (resp. $a_n\to y$), and with $a_n<x_n$ for all $n\ge1$.
It is plain that $M^{x_n}_t\to M^x_t$ and $M^{a_n}_t\to M^a_t$ a.s. and in $L^1$. Moreover for each $n\ge1$,
\begin{align*}
\E[M^{x_n}_t-M^{a_n}_t|\G^x]=\E\left\{\E[M^{x_n}_t-M^{a_n}_t|\G^{x_n}]|\G^x\right\}
=0,
\end{align*}
hence taking the limit as $n\to\infty$ in that identity, we deduce that \eqref{eq:mart} holds true for any $0<a<x$.
The result is established.
\epf

We now show that  for any $A\in\EE$, the Borel field of $E$, any $t>0$, 
\begin{equation}\label{Qchar}
\int_A\mathbf Q_{0,t}(du)=\int_Au(t)Q(du),
\end{equation}
where the $\sigma$--finite measure $Q$ on $(E,\EE)$  is the excursion measure of Feller's critical diffusion \eqref{CFD}, in the sense of Pitman and Yor, see formula (3a) in \cite{PY}, with the scale function $s$ being chosen as $s(y)=y$.
To see \eqref{Qchar}, first note that for all $\Phi\in C_b(E)$,
 \begin{align*}
\lim_{y \to 0}\frac{1}{y}\E_y[\Phi(U)]=\int_E\Phi(u)Q(du).
\end{align*}
 Since $\frac{1}{y}\E_y[U_t] = 1$ for all $y >0$, this implies by a  uniform integrability argument  that
 \begin{align*}
\lim_{y \to 0}\frac{1}{y}\E_y[\Phi(U)U_t]=\int_E\Phi(u)u(t)Q(du).
\end{align*}
By the definition of $\mathbf Q_{y,t}$ the l.h.s is $\lim_{y \to 0}\int_E\Phi(u)\mathbf Q_{y,t}(du)$, which is $\int_E\Phi(u)\mathbf Q_{0,t}(du)$ in view of \eqref{Fimm}. This proves \eqref{Qchar}.

  The measures $Q\circ u_t^{-1}$, $t > 0$, constitute an entrance law of  of the Feller diffusion \eqref{CFD}. This entrance law, which also figures in  formula (3.2) of \cite{PY},  is  given by 
\begin{equation}\label{Qent}
Q\circ u_t^{-1}=(2t)^{-1}\text{Exp}((2t)^{-1}), t > 0.
\end{equation}
 Indeed, it is readily checked from formula \eqref{FDimmt} that the distribution of $u_t$ under $\mathbf Q_{y,t}$ is Gamma$(2,2t)$, which is the size-biasing of $\text{Exp}((2t)^{-1})$. On the other hand, it is immediate from \eqref{Qchar} that $Q_{0,t}\circ u_t$ is the size-biasing of $Q\circ u_t^{-1}$. From this the claim \eqref{Qent} is immediate. Let us also note that our probabilities $\mathbf Q_{y,\infty}$ are the ``upward diffusions'' $P_y^{\uparrow}$ of~\cite{PY}.
 
Combining \eqref{Mx}, \eqref{Qchar} and Remark \ref{pathidentity}, we immediately arrive at our main result
\begin{theorem}\label{th:path-val-SDE}
The path--valued process $\{Z^x_\cdot,\ x>0\}$ admits the decomposition
\begin{equation}\label{eq:E}
Z^x =\int_{[0,x]\times E}u L(Z^\xi,u) Q(du)d\xi+M^x,
\end{equation}
where $M^x$ is a $C([0,+\infty);\R)$--valued c\`adl\`ag martingale (if $C([0,+\infty);\R)$ is equipped with the topology of uniform convergence on compacts). 
\end{theorem}

We know that $x\mapsto Z^x$ arises as a sum of excursions, as was stated above in Corollary~\ref{corRC}. Call $N_Z(d\xi,du)$ the corresponding point process, which is such that for all $x>0$,
$$Z^x=\int_{[0,x]\times E}uN_Z(d\xi,du).$$
The above statement shows that the predictable intensity measure of $N_Z$ is
$$L(Z^\xi, u)Q(du)d\xi.$$
Intuitively (and somewhat informally stated) this means that, given $(Z^\xi)_{0\le \xi < x}$,   the predicted increment of $Z$ in the next bit $dx$ of ancestral mass is a Poisson point process with intensity measure $L(Z^{x-},u)Q(du)\, dx $. This is made precise by the following statement, which was conjectured in the case of a logistic drift in \cite{PW1}:
\begin{corollary}\label{cor:martpb}
For bounded $g: \mathbb R_+ \to \mathbb R_+$ and $z\in E$, put $\Phi_g(z) := e^{-\langle g, z\rangle}$. Then, for this class of functions,
$$A\Phi_g(z) := \Phi_g(z) \int_E\left(e^{-\langle g, u\rangle}-1\right)  L(z,u) Q(du)$$
gives the generator of $Z$ in the sense that for all  $g: \mathbb R_+ \to \mathbb R_+$, 
\begin{equation}\label{mart}
\Phi_g(Z^x)- \Phi_g(Z^0)-\int_{[0,x]\times E} A\Phi_g(Z^\xi) d\xi, \, x \ge 0 \qquad \mbox{ is a martingale. }
\end{equation}
\end{corollary}
\bpf
The validity of \eqref{mart} can be seen by writing
\begin{align*}
\Phi_g(Z^x)- \Phi_g(Z^0) &= \int_{[0,x]\times E} 
(\Phi_g(Z^{\xi-}+u)-\Phi_g(Z^{\xi-}))
N(d\xi, du)
\\ & = \int_{[0,x]\times E}
e^{-\langle g, Z^{\xi-}\rangle}
(e^{-\langle g, u \rangle}-1) N(d\xi, du).
\end{align*}
The same line of arguments that led to \eqref{eq:E} now shows that the r.h.s. equals
\begin{align*}
\int_{[0,x]\times E} 
e^{-\langle g, Z^{\xi-}\rangle}
(e^{-\langle g, u \rangle}-1)
L(Z^{\xi-},u)Q(du) d\xi + M^x_g, \, x \ge 0,
 \end{align*}
 for some real-valued c\`adlag martingale $\{M_g^x,\ x>0\}$, which yields \eqref{mart}.
 \epf
 \\
{\bf Acknowledgement}

\smallskip
 \noindent We thank a referee for valuable comments that led to an improvement of the paper.

\end{document}